\documentclass[11pt]{article}
\input epsf
\usepackage{amsmath}
\usepackage{amssymb}
\usepackage{amscd}
\usepackage{epsfig}
\usepackage{psfrag,color}

\newtheorem{dfn}{Definition}[section]

\newtheorem{cor}[dfn]{Corollary}

\newtheorem{lem}[dfn]{Lemma}
\newtheorem{prop}[dfn]{Proposition}

\newtheorem{rem}[dfn]{Remark}

\newtheorem{thm}[dfn]{Theorem}

\def\lra{\longrightarrow}
\def\lmt{\longmapsto}

\def\ra{\rightarrow}

\def\C{{\mathbf C}}
\def\R{{\mathbf R}}

\def\3{\ss}
\def\8{\infty}

\def\qed{$\square$\bigskip}

\newcommand{\Isom}{\operatorname{Isom^+}}
\newcommand{\isom}{\mathfrak{isom^+}}

\newcommand{\so}{\mathfrak{so}}
\newcommand{\E}{\mathcal{E}}
\newcommand{\Fold}{\mathcal{F}}

\begin{document}

\title{Deforming Euclidean cone 3-manifolds}
\author{Joan Porti\footnote{Partially supported by the 
Spanish MEC through grant BFM2003-03458 and European funds FEDER} and Hartmut Weiss
\footnote{Partially supported by DFG Schwerpunktprogramm 1154 Globale Differentialgeometrie}
}
\date{\today}

\maketitle

\begin{abstract}
Given a closed orientable Euclidean cone 3-manifold $C$ with cone angles $\leq
\pi$ and which is not
almost product, we describe the space of constant curvature cone structures on $
C$ with cone angles $<\pi$.
We establish a regeneration result for such Euclidean cone manifolds into
spherical or hyperbolic ones and we also deduce global rigidity for Euclidean
 cone structures.
\end{abstract}

%\tableofcontents

\section{Introduction}
Let $C$ be a closed orientable Euclidean cone 3-manifold with cone angles
$\leq\pi$. Its singular locus $Sing(C)$
is a trivalent graph consisting of $q$ circles and edges, that we enumerate from
 $i=1$ to $q$. To those edges one can associate a
multiangle $\bar \alpha=(\alpha_1,\ldots,\alpha_q)$, where $\alpha_i$ is the
angle of the $i$-th circle or edge of $Sing(C)$. Let further $\bar l =(l_1, \ldots,
 l_q)$ denote the vector of lengths of singular circles and edges. We are
interested in the space of cone structures on $C$ obtained by changing the
multiangle.

\begin{dfn}
 A \emph{hyperbolic cone structure} on $C$ is a hyperbolic cone
 manifold $X$ with an embedding $i\! :\! (X,Sing(X))\to
 (C,Sing(C))$ such that $i(X^{Smooth})$ is a retract of
 $C^{Smooth}$.
\end{dfn}

\begin{dfn}
A Euclidean cone 3-manifold $C$ is \emph{almost product} if
$C=(S^1\times E^2)/G$ where $E^2$ is a  two-dimensional cone
manifold and $G<\operatorname{Isom} (S^1)\times
\operatorname{Isom}(E^2)$ is a finite subgroup.
\end{dfn}

\begin{thm}
\label{thm:main}
Let $C$ be a closed orientable Euclidean cone 3-manifold with
cone angles $\leq    \pi$. If $C$ is not almost product, then for
every multiangle $\bar\alpha\in (0,\pi)^{q}$ there exists a unique
cone manifold structure of constant curvature in $\{-1,0,1\}$ on
$C$ with those cone angles.

If all cone angles of $C$ are $\pi$, then every point in $(0,\pi)^q$ is the multiangle
of a hyperbolic cone structure on $C$.

If some of the cone angles is $<\pi$, then
the subset $E\subseteq (0,\pi)^q$   of
multiangles of Euclidean cone structures is a smooth, properly embedded hypersurface
that splits $(0,\pi)^{q}$ into 2 connected components $S$ and $H$,
corresponding to multiangles of spherical and hyperbolic cone
structures respectively. Furthermore, for each $\bar \alpha \in E$ the tangent space of $E$ at $\bar \alpha$ is orthogonal to the vector of singular lengths $\bar{l}$.
\end{thm}

\noindent This theorem describes the structures corresponding to multiangles in
the open cube. We can also describe the structures corresponding to multiangles
contained in some parts of the boundary. For instance, we show that the
multiangles contained in $\partial [0,\pi)^q\cap
[0,\pi)^q$ (i.e.~none of the angles is $\pi$ and at least one of them is
zero) are angles of hyperbolic cone structures. However, we do not give a description
of structures in all points of $\partial [0,\pi]^q$, as this would involve
studying collapses at cone angle $\pi$ and deformations and regenerations of
Seifert fibered geometries.

% We can describe what happens in some subsets of the boundary
%$\partial[0,\pi]^q$. For instance,
%Using Proposition~\ref{prop:largerspherical}, we can further prove:
% the multiangles
%in the closure of the spherical region $S$ in
%$(0,\pi]^q-\overline E$ are cone angles of a spherical cone
%manifold. Thus we get:

\begin{cor}
\label{cor:sphericity}
 Let $\mathcal O$ be an irreducible closed
orientable $3$-orbifold. If there exists a Euclidean cone
structure $C$ on $\mathcal O$ with cone angles strictly less than
the orbifold angles of $\mathcal O$, then $\mathcal O$ is
spherical.
\end{cor}

\begin{proof}
To prove this corollary, we  show in Lemma~\ref{lem:Cnotap} that
the Euclidean cone structure is not almost product, using the
irreducibility of $\mathcal O$.
By Theorem~\ref{thm:main} we can obtain a spherical cone structure on $C$ by increasing any of the cone angles. We are using here that $E$ is orthogonal to the
vector of singular lengths $\bar l$. By Proposition~\ref{prop:largerspherical}
the orbifold angles can be realized by a spherical cone structure.\hfill\qed
\end{proof}

\noindent This corollary gives an alternative argument to the last step in the
proof
of the orbifold theorem in \cite{BLP2}, which is more natural from the point
of view of cone manifolds. D.~Cooper and S.P.~Kerckhoff have announced a
different approach to the spherical uniformisation.

%
% except in
%the case of Proposition~\ref{prop:largerspherical}, giving
%spherical structures.
%

%We shall prove in Proposition \ref{prop:lnormal} that the codimension
%one submanifold $E\subset (0,\pi)^{q}$ of multiangles of Euclidean
%cone structures has the following property: if $\bar\alpha\in E$
%and $\bar l=(l_1,\ldots,l_q)$ are the corresponding lengths of
%edges and circles of $Sing(C)$, then $E$ is orthogonal to $\bar l$
%at $\bar \alpha$.

A special case of this theorem was proved in \cite{Porti}, assuming
that the singular locus was a knot and a technical hypothesis involving
cohomology. Even if some of the techniques of \cite{Porti} are used,
this paper does not rely on it.

Now we describe the organization of the paper: We are interested in
the rotational part of the holonomy of the Euclidean cone manifold
in $SO(3)$ (in fact its lift to $Spin(3)$), that we denote by $\rho_0$.
Some basic properties of this representation are studied in
Section~\ref{sec:almostproduct}. In order to deform the structure,
we view $SO(3)$ as the stabilizer of a point in $\mathbf H^3$ or
$\mathbf S^3$, thus we study the local properties of the varieties
of representations in the isometries of $\mathbf H^3$ and $\mathbf S^3$ around $
\rho_0$.
 This is done
in  Section~\ref{sec:varofrep}, using the cohomological results of
Section~\ref{sec:cohomologytangentbundle}.
In Section~\ref{section:regenerationstructures} we give the conditions
for a deformation of the representation $\rho_0$
to correspond to a path of hyperbolic or spherical cone manifolds. This is
applied in
Section~\ref{sec:pathsofstructures} to construct paths of hyperbolic and
spherical
structures by deforming one of the
cone angles. The analysis of the local parametrization of
the variety of representations is completed in
Section~\ref{sec:foldlocus}, where all deformations of the structure are
constructed. The global results are established in Section~\ref{sec:global}. Finally,
in Section~\ref{sec:example} we illustrate the main theorem for cone manifold
structures on the 3-sphere with singular locus given by the Whitehead link.

\paragraph*{Acknowledgments} We thank Bernhard Leeb for helpful conversations
on the contents of this paper and for sharing his ideas with us.
%The authors are partially supported by Spanish MCYT through grant
%BFM2003-03458 and by DFG Schwerpunktprogramm 1154 Globale Differentialgeometrie.

\section{Almost product structures}
\label{sec:almostproduct}

Given a Euclidean cone manifold $C$, its smooth part is denoted by
$M=C-Sing(C)$. Then $M$ has a non-complete and non-singular flat
Riemannian metric. In particular it has a developing map
\[
dev\!:\widetilde M\to \mathbf R^3
\]
and a holonomy representation
\[
hol\!:\pi_1M\to\Isom\R^3=SO(3)\ltimes \R^3\,.
\]
We consider the composition of $hol$ with the projection \(ROT \!:
\Isom\R^3 \to SO(3)\).

\begin{lem}
\label{lem:almostproduct}
 $C$ is almost product if and only if the image of
$ROT\circ hol\! :\! \pi_1M\to SO(3) $ is contained in $O(2)\subset
SO(3)$.
\end{lem}

\begin{proof}
If  $C$ is almost product then the lemma is straightforward. To prove the
converse, we assume that the image of $ROT\circ hol$ is contained in
$O(2)\subset SO(3)$. Since $O(2)$ has two components, composition with
$O(2)\to\pi_0 O(2)$ defines a morphism $\pi_1M\to\mathbf  Z/2\mathbf  Z$
and, if it is non-trivial, we take the corresponding covering. This
induces a branched covering of $C$, with branching locus contained in the
edges and circles of $Sing(C)$ that have cone angle $\pi$. Hence we may
assume that the image of $ROT\circ hol$ is contained in $SO(2)$. The
vertical foliation of $\mathbf  R^3$ invariant under $SO(2)$ induces an
isometric flow of $M$.  It extends to a flow on the whole of $C$, since this
restriction of the holonomy implies that the components of $Sing(C)$ are
parallel to the leaves of the flow.

If $Sing(C)=\emptyset$, then  $C$ is a smooth flat manifold and
the lemma is well known. Hence we assume that
$Sing(C)\neq\emptyset$. Since the components of $Sing(C)$ are
leaves of the flow, either all leaves are closed or they
accumulate on tori at constant distance from the components of
$Sing(C)$. The latter is not possible, since the concentric tori
of accumulation of fibers have principal curvatures of constant
sign. Hence the leaves of the flow are closed, and they are leaves
of a Seifert fibration, whose basis has a natural structure of a Euclidean cone
2-manifold.

One can construct a horizontal surface $H$ perpendicular to the
flow, and we can enlarge it as much as possible to be complete.
Since $H$ has cone angles $\leq \pi$, the singular locus of $H$ is
finite. Thus $H$ is compact and therefore the Seifert fibration is
virtually a product. \qed
\end{proof}

\begin{dfn}
The Euclidean cone manifold $C$ is \emph{Seifert fibered} if it
admits a Seifert fibration such that $Sing(C)$ is a union of
fibers. It is \emph{almost Seifert fibered} if it is the quotient
of a Seifert fibered cone manifold by a finite subgroup of
isometries that preserve the fibration.
\end{dfn}

\noindent A cone manifold which is almost Seifert fibered but not Seifert
fibered admits a partition by circles
and intervals, so that the end-points of the intervals lie in
singular edges with cone angle $\pi$. This corresponds to an orbifold Seifert
fibration. In particular, if cone angles are $<\pi$, then almost Seifert
fibered implies Seifert fibered.

\begin{cor}\label{cor:almostseifert}
The Euclidean cone manifold $C$ is almost product iff it is almost
Seifert fibered. In particular, if all cone angles of $C$ are $<\pi$
and its smooth part is not Seifert fibered, then $C$ is not almost product.
\end{cor}

\begin{proof}
Almost product implies trivially almost Seifert fibered. Assume
now that $C$ is Seifert fibered. Then the smooth part of $C$ is a
product. Let $f\in\pi_1(C-Sing(C))$ represent the fiber of this
fibration. Since the components of $Sing(C)$ are also fibers, the
holonomy of $f$ is either a nontrivial translation or a screw
motion with non-trivial translation length. Since $f$ is central,
the direction of this vector must be preserved. Thus the image of
$ROT\circ hol$ is contained in $SO(2)$. It follows that when $C$
is almost Seifert fibered, then the image of $ROT\circ hol$ is
contained in $O(2)$, and by Lemma~\ref{lem:almostproduct} $C$ is
almost product. \hfill\qed
\end{proof}

\noindent We consider the action by conjugation of $SO(3)$ on $M_{3\times 3}(\mathbf
R)=End(\mathbf  R^3)$, the vector space of $3\times 3$ matrices with real
entries.

\begin{lem}
\label{lem:invariantsubspace}
 $C$ is not almost product if and only if  the only subspace
of $M_{3\times 3}(\mathbf  R)$ invariant under $
 ROT\circ hol$
is the space of diagonal matrices.
\end{lem}

\begin{proof}
As $SO(3)$-module acting by conjugation, we have the splitting:
\[
M_{3\times 3}(\mathbf  R)=\mathbf  R\oplus so(3)\oplus Z^5,
\]
where $\mathbf  R$ is the set of diagonal matrices, $so(3)\cong\mathbf
R^3$ is the subspace of antisymmetric ones and $Z^5\cong \mathbf  R^5$ is
the subspace of traceless symmetric matrices.

The action on $\mathbf  R$ is trivial and this space is always invariant.
The action on $so(3)$ is equivalent to the usual action of $SO(3)$ on
$\mathbf  R^3$, and having a nonzero invariant vector in $\mathbf  R^3$ is
equivalent for $ROT\circ hol$ to be contained in $SO(2)$.

 Since symmetric matrices
diagonalize orthogonally, matrices in $Z^5$ have an orthogonal
basis of eigenvectors. Furthermore the sum of eigenvalues is
zero, and therefore every non-trivial element has at  least a one
dimensional eigenspace.
 Thus having a non-trivial invariant element of $Z^5$ is equivalent for
$ROT\circ hol$ to have an invariant line, i.e.~to be
contained in $O(2)$. \hfill\qed
\end{proof}

\begin{lem}
 \label{lem:Cnotap}
 Let $C$ and $\mathcal O$ be as in the situation of
Corollary~\ref{cor:sphericity}. Then $C$ is not almost product.
\end{lem}

\begin{proof} By hypothesis the cone angles of $C$ are $<\pi$.
Therefore if $C$ is almost product, then $C$ is a quotient of
$E^2\times S^1$ where $E^2\cong S^2(\alpha,\beta,\gamma)$ is a
Euclidean turnover (a cone manifold structure on $S^2$ with three
cone points of cone angles satisfying $\alpha+\beta+\gamma=2\pi$).
This cone manifold $E^2$ embeds in $C$ and induces an essential
spherical 2-orbifold in $\mathcal O$, contradicting the
irreducibility of $\mathcal O$. \hfill\qed
\end{proof}

\section{Cohomology of the tangent bundle}
\label{sec:cohomologytangentbundle}
We assume from now on that $M$ is the smooth part of a closed
Euclidean cone manifold $C$ which is not almost product. The
latter assumption may be viewed as a nondegeneracy condition and
is in that respect similar to the assumption of not being Seifert
fibered in the deformation theory of spherical cone-manifolds,
cf.~\cite{BLP2}, \cite{Wei1}.

In the Euclidean case the flat bundle of infinitesimal isometries
$\E=\mathfrak{so}(TM)\oplus TM$ contains the bundle of
infinitesimal translations $\E_{trans}=\{0\} \oplus TM$ as a
parallel subbundle. The inclusion $\E_{trans} \subset \E$ fits
into a short exact sequence of flat vectorbundles and
connection-preserving maps
$$
0 \ra \E_{trans} \ra \E \ra \E_{rot} \ra 0 \,,
$$
where $\E_{rot}=\mathfrak{so}(TM)\oplus\{0\}$ carries the {\em
quotient} connection ($\E_{rot} \subset \E$ is {\em not}
parallel!). In fact $\E_{rot}$ is again isomorphic to $TM$ as a
flat bundle. This can for instance be seen if we look at the
corresponding short exact sequence of $\Isom \R^3$-modules
$$
0 \lra \R^3_{Ad(\,\cdot\,)|_{\R^3}} \lra \isom \R^3_{Ad} \overset{ROT_*}{\lra} \
so(3)_{Ad \circ ROT} \lra 0\,.
$$
Now we have isomorphisms of $\Isom \R^3$-modules $\so(3)_{Ad \circ
ROT} \cong \R^3_{Ad(\,\cdot\,)|_{\R^3}} \cong \R^3_{ROT}$, which
in turn yield isomorphisms of flat bundles $\E_{rot} \cong
\E_{trans} \cong TM$. Note that the first isomorphism is
particular to our 3-dimensional situation. We will freely identify
$TM$ with either $\E_{trans}$
or $\E_{rot}$ from now on.\\

\noindent The translational part of a Euclidean isometry $\phi$
w.r.t.~a base point $p \in \R^3$
$$
trans_p(\phi)=\phi(p)-p
$$
gives rise to a cocycle
$
trans_p : \Isom \R^3 \ra \R^3
$
twisted by $ROT$, i.e.~
$$
trans_p(\phi_1\phi_2) = trans_p(\phi_1) + ROT(\phi_1)(trans_p(\phi_2))
$$
for $\phi_1,\phi_2 \in \Isom \R^3$.

\begin{dfn}
The \emph{holonomy cocycle} is the composition
$$trans_p\circ hol: \pi_1M \lra \R^3$$
for some fixed $p\in\mathbf  R^3$.
\end{dfn}

\noindent This is a cocycle twisted by $ROT\circ hol$, i.e.~
\[
trans_p\circ hol(\gamma_1\gamma_2)= trans_p\circ hol(\gamma_1)+
(ROT\circ hol(\gamma_1))(trans_p\circ hol(\gamma_2)).
\]
for $\gamma_1, \gamma_2 \in \pi_1M$. The cohomology class of this
cocycle is independent of $p$ and is denoted by
\[
v=[trans_p\circ hol]\in H^1(\pi_1M; \R^3_{ROT \circ hol})\,.
\]
Let $\omega \in \Omega^1(M;TM)$ be closed and $\gamma$ a loop
based at $x \in M$. We define the integral of $\omega$ along
$\gamma$ as
$$
\int_{\gamma} \omega = \int_0^1 \tau^{-1}_{\gamma(t)} \omega(\dot{\gamma}(t)) dt
 \in \R^3 \,,
$$
where $\tau_{\gamma(t)}$ denotes parallel transport along $\gamma$
from $x=\gamma(0)$ to $\gamma(t)$ and we identify $T_xM$ with
$\R^3$ by developing $M$ on $\R^3$. For $\omega \in
\Omega^1(M;TM)$ closed the assignment $\gamma \mapsto
\int_{\gamma}\omega$ defines a group cocycle $z_{\omega} \in
Z^1(\pi_1M;\R^3_{ROT \circ hol})$,  the so-called period cocycle
of $\omega$. The period map $\omega \mapsto z_{\omega}$ descends
to an isomorphism in cohomology
\begin{align*}
P : H^1(M;TM) &\overset{\cong}{\lra} H^1(\pi_1 M; \R^3_{ROT \circ hol})\\
[w] & \lmt [z_{\omega}]\,.
\end{align*}
Observe that on each side there is a canonically defined cohomology class. On
 the left hand side we have $id \in \Omega^1(M;TM) = \Gamma(M;T^*M \otimes
 TM)$. $id$ is parallel and gives rise to a class $[id] \in H^1(M;TM)$. On the
 right hand side we have the class of the holonomy cocycle $v$. In fact, these
 classes coincide modulo the period map and correspond geometrically to a global
 rescaling of $M$:

\begin{lem}
\label{lem:Pid=v}
$P([id])=v$ and the image of $v$ in $H^1(\pi_1M;\isom \R^3_{Ad \circ hol})$ is
  tangent to a deformation, which globally rescales the Euclidean metric on $M$.

\end{lem}

\begin{proof}
Let $\gamma$ be a loop based at $x \in M$ and $\tilde{\gamma}$ a lift of
$\gamma$ to $\widetilde{M}$. We compute
$$
\int_{\gamma} id = trans_{dev(\tilde{\gamma}(0))} \circ hol(\gamma)\,,
$$
which shows that the periods of $id$ reproduce the holonomy cocycle defined with

respect to $dev(\tilde{\gamma}(0)) \in \R^3$.
Furthermore for $\gamma \in \pi_1M$ let
$$
\rho_t(\gamma) = (ROT \circ hol (\gamma), (t+1) trans_0 \circ hol(\gamma)) \in
\Isom \R^3\,.
$$
Then $\rho_0 = hol$ and $\rho_t$ is the  holonomy representation
of $M$ with rescaled Euclidean metric $(t+1)^2g$. We obtain $
\rho_t(\gamma)hol(\gamma)^{-1} = ( 1 , t\cdot trans_0 \circ
hol(\gamma)) $ and consequently $
\left.\textstyle\frac{d}{dt}\right\vert_{t=0}
\rho_t(\gamma)hol(\gamma)^{-1} = (0, trans_0 \circ hol(\gamma)) $.
This means the corresponding path of characters is precisely
tangent to $v$, cf.~the later discussion in Section
\ref{sec:varofrep}. \hfill\qed
\end{proof}

\noindent Let $U_{\varepsilon}(Sing(C))$ denote the $\varepsilon$-tube around
$Sing(C)$ in $C$ minus the singular locus itself. Let $\overline{M}$
denote a compact core of $M$. $\partial\overline{M}$ is a disjoint union of
tori and surfaces of higher genus. Note that the pair $(M,U_{\varepsilon}(Sing(C
)))$
is homotopy equivalent to the pair $(\overline{M},\partial\overline{M})$, by
homotopy invariance we may identify the corresponding cohomology groups.\\

\noindent If $(r_i,\theta_i,z_i)$ are cylindrical coordinates
around  the $i$-th edge and $\varphi_i=\varphi_i(z_i)$ is a
function with $\varphi_i\equiv 0$ near $z_i=0$ and
$\varphi_i\equiv l_i$ near $z_i=l_i$, we define as in \cite{Wei1}
$$
\omega^i_{len} = d\varphi_i \otimes \textstyle
\frac{\partial}{\partial z_i} \in \Omega^1(U_{\varepsilon}(Sing(C));TM)\,.
$$
As in \cite{Wei1} we obtain:

\begin{lem}\label{L2}
The forms $\omega^i_{len}$ are bounded on $U_{\varepsilon}(Sing(C))$,
 in particular $L^2$.
\end{lem}

\noindent Geometrically, the classes $[\omega^i_{len}]$ correspond to
deformations of $U_{\varepsilon}(Sing(C))$ changing the length of the $i$-th
edge. These deformations are independent. Therefore the following does not
come as a surprise:

\begin{lem}\label{lem:nontriv}
The classes $[\omega^i_{len}]$ are linearly independent in
$H^1(\partial \overline{M};TM)$, furthermore the relation
$[id]|_{\partial \overline{M}}=\sum_i [\omega^i_{len}]$ holds.
\end{lem}

\begin{proof}
The first assertion follows as in \cite{Wei1}.
The periods of both forms $id$ and $\sum_i \omega^i_{len}$ reproduce the
holonomy cocycle restricted to $\pi_1 \partial \overline{M}$. Since the period
map descends to an isomorphism in cohomology, the result follows. \hfill \qed
\end{proof}

\begin{lem} $H^0(M;TM)=0$.
\end{lem}

\begin{proof} This group is isomorphic to the subspace of $\mathbf
R^3$ invariant under $ROT\circ hol$, which is trivial since we assume
that $C$ is not almost product (see the proof of
Lemma~\ref{lem:invariantsubspace}). \hfill\qed
\end{proof}

\noindent From the long exact cohomology sequence we obtain:

\begin{cor}
$H^0(M;\E)=0$ and the map $H^1(M;TM) \ra H^1(M;\E)$ is an injection.
\end{cor}

\noindent To compute $H^1(M;TM)$ we use the $L^2$-cohomology:
$H^1_{L^2}(M;TM)$. Our starting point is the following theorem of
\cite{Wei1}, where the hypothesis about cone angles is used:

\begin{thm}
    \label{thm:L2}
If $M$ is the smooth part of a closed Euclidean cone manifold with
cone angles $\leq \pi$, then
$$
H^1_{L^2}(M;TM) \cong \{ \omega \in \Omega^1(M;TM) \vert \nabla \omega = 0 \}\,.
$$
\end{thm}

\begin{cor}
\label{cor:l2}
 $H^1_{L^2}(M;TM)\cong\mathbf  R$.
\end{cor}

\begin{proof} Since every element of $H^1_{L^2}(M;TM)$ is
represented by a parallel form, this cohomology group is isomorphic to the
space of elements in $(\mathbf  R^3)^*\otimes \mathbf  R^3=End(\mathbf
R^*)\cong M_{3\times 3}(\mathbf  R^3)$ which are invariant under the action
of $ROT\circ hol$ by conjugation. By Lemma~\ref{lem:invariantsubspace},
this subspace is one dimensional. \hfill\qed
\end{proof}

\noindent Let $\overline M$  denote a compact core of $M$. By homotopy
invariance $H^*(M;TM)\cong H^*(\overline M;T M)$.

\begin{lem}
\label{lem:inj}
 The map \(H^1(M; T M)\to H^1(\partial \overline M;
T M)\) is an injection.
\end{lem}

\begin{proof}
Looking at the exact sequence of the pair $(\overline M,\partial
\overline M )$, it suffices to show that the image of the map $j\!
: \! H^1(\overline M,\partial \overline M;TM)\to H^1(\overline
M;TM)$ is zero. To prove it, we observe that $j$ factors through
$H^1_{L^2}(M;TM)$. By Corollary~\ref{cor:l2}, $H^1_{L^2}(M;TM)$ is
one-dimensional and its generator is represented by the identity
of $TM$. Therefore the image of $j$ is contained in the span of
$[id]$. By Lemma~\ref{lem:nontriv}, $[id]$ restricts to a
nontrivial class in $H^1(\partial \overline{M};TM)$. On the other
hand, by exactness the image of $j$ is contained in the kernel of
$H^1(M;TM) \ra H^1(\partial\overline{M};TM)$, hence it is zero.
\hfill \qed
\end{proof}

\noindent As a corollary of the proof above we notice:

\begin{cor}\label{L2:inj}
The map $H^1_{L^2}(M;TM) \ra H^1(M;TM)$ is an injection.
\end{cor}

\noindent We may thus view $L^2$-cohomology as a subspace of ordinary cohomology
 in
degree 1.

\begin{lem}
    \label{lem:H1q}
$\dim H^1(M; T M)=q$.
\end{lem}

\begin{proof}
By Lemma~\ref{lem:inj}, the long exact sequence of the pair and Poincar\'e
duality, we have a short exact sequence:
\[
0\to H^1(M; T M)\to H^1(\partial \overline M; T M)\to
H^2(\overline M,\partial \overline M;TM)\to 0
\]
Since $H^1(M; T M)$ and $H^2(\overline M,\partial \overline M;TM)$
are Poincar\'e dual,
\[
\dim H^1(M; T M)=\frac12\dim  H^1(\partial \overline M; T M)\,.
\]
In addition, $\dim  H^1(\partial \overline M; T M)=2 q$ (see
\cite{BLP2} or \cite{Wei1}). \hfill\qed
\end{proof}

\noindent Continuing the long exact sequence we get:

\begin{cor}
    \label{lem:H2isomorphic}
The inclusion induces an isomorphism $H^2(M;TM)\cong
H^2(\partial\overline M;TM)$.
\end{cor}

\noindent Let $H^*(\pi_1M;\R^3_{ROT\circ hol})$ denote the cohomology of the
$\pi_1M$-module $\mathbf  R^3$ with the action induced by $ROT\circ hol$.
Notice that the induced flat bundle over $M$ is isomorphic to $TM$.

\begin{lem}
 \label{lem:cohomologies}
 There is a natural map
$$
H^i(\pi_1M; \R^3_{ROT\circ hol}) \to H^i(M;TM)
$$
which is an isomorphism for $i=0,1$ and an injection for $i=2$.
\end{lem}

\begin{rem}
When $i=1$, this map is the inverse of the isomorphism $P$ in Lemma~\ref{lem:Pid=v}.
\end{rem}

\begin{proof}
There is a natural isomorphism  $H^i(\pi_1M;\R^3_{ROT\circ
hol})\cong H^i(K(\pi_1M,1);\R^3_{ROT\circ hol})$, where
$K(\pi_1M,1)$ denotes any aspherical CW-complex with the same
fundamental group as $M$, and $\R^3_{ROT\circ hol}$ also denotes
the flat bundle over $K(\pi_1M,1)$ with fiber $\mathbf R^3$ and
holonomy $ROT\circ hol$. The space $K(\pi_1M,1)$ can be
constructed from $M$ itself by adding cells of dimension $\geq 3$.
Thus the cohomology of the pair $(K(\pi_1M,1),M)$ is trivial in
dimension $\leq 2$, and the lemma follows easily from the long
exact sequence in cohomology for the pair $(K(\pi_1M,1),M)$.
\hfill \qed
\end{proof}

\begin{cor}
   \label{cor:injectiontoboundary}
    Let $\partial_1\overline M,\ldots,\partial_k\overline M$ denote the components
    of $\partial \overline M$. The restriction induces an injection:
    $$
    0\to H^2(\pi_1M; \R^3_{ROT\circ hol})\to \bigoplus_{i=1}^k H^2(\pi_1\partial
_i\overline M; \R^3_{ROT\circ hol})
    $$
\end{cor}

\begin{proof}
The result follows from Corollary~\ref{lem:H2isomorphic}, Lemma~\ref{lem:cohomologies}
and the asphericity of the surfaces $\partial_i\overline M$.
\hfill\qed
\end{proof}

\section{Varieties of representations}\label{sec:varofrep}
Since $M$ is three dimensional by a theorem of Culler, cf.~\cite{Cul}, $hol$
lifts to a representation
$$
\rho=\widetilde{hol}: \pi_1M \rightarrow \widetilde{\Isom \R^3} \cong SU(2) \ltimes \R^3\,,
$$
and similarly $ROT\circ hol$ lifts to
\[
\rho_0=\widetilde{ROT}\circ \widetilde{hol}\!:\! \pi_1M\to \widetilde {SO(3)}\cong SU(2)\,,
\]
where $\widetilde{ROT}: \widetilde{\Isom}\R^3 \rightarrow SU(2)$ is a lift of
$ROT$. We will drop the notational distinction between $\widetilde{hol}$ and
$hol$, resp.~between $\widetilde{ROT}$ and $ROT$ from now on.

Let $R(M, SU(2))=R(\pi_1M,SU(2))$ and $R(M, SL_2(\mathbf  C))=R(\pi_1M, SL_2(\mathbf  C))$
denote the varieties of
representations in $SU(2)$ and $SL_2(\mathbf  C)$ respectively. Let $X(M,
SU(2))$ and $X(M,SL_2(\C))$ be the varieties of characters. Let
$\chi_0\in X(M, SU(2))\subset X(M,SL_2(\mathbf  C))$ denote
the character of $\rho_0$.

We recall some general facts about spaces of representations:

\noindent Let $\Gamma$ be a finitely generated discrete group and
$G$ a Lie group. We equip $R(\Gamma, G)$ with the compact-open
topology. Once and for all we fix a presentation $\langle
\gamma_1,\ldots,\gamma_n | (r_i)_{i\in I}\rangle$. Without loss of
generality we may assume that the set of generators contains with
any of its elements its inverse. The relations $r_i$ determine
functions $f_i: G^n \rightarrow G$ which identify $R(\Gamma,G)$
with the subset $\{ f_i = 1 \} \subset G^n$.

If $G$ is a real algebraic group (as for example $G=SU(2)$ or
$\widetilde{\Isom \R^3}$), the $f_i$ are polynomial functions and
$R(\Gamma,G)$ acquires the structure of a real algebraic set.
Similarly, if $G$ is a complex algebraic group (for example
$G=SL_2(\C)$), then $R(\Gamma,G)$ is a complex algebraic set.

Let $G\ltimes_{Ad} \mathfrak{g}$ be the semidirect product of $G$ and $\mathfrak
 g$
defined via $Ad: G \rightarrow Aut(\mathfrak{g})$. Let $\rho \in R(\Gamma, G)$
be given. A simple but important observation is that $z \in
Z^1(\Gamma;\mathfrak{g}_{Ad \circ \rho})$ if and only if $(\rho,z) \in
R(\Gamma, G \ltimes_{Ad} \mathfrak{g})$. The relations therefore determine
functions $g_i: G^n \times \mathfrak{g}^n \rightarrow \mathfrak{g}$ such
that $R(\Gamma, G \ltimes_{Ad} \mathfrak{g}) = \{ (f_i,g_i) = (1,0) \}$ and
$Z^1(\Gamma; \mathfrak{g}_{Ad \circ \rho}) = \{ g_i(\rho, \cdot) = 0 \}$. Note
that for $\rho$ fixed the functions $g_i(\rho, \cdot): \mathfrak{g}^n
\rightarrow \mathfrak{g}$ are linear. If we identify $T_AG$ with
$\mathfrak{g}$ via right translation, then we have:

\begin{lem}
$g_i(\rho, \cdot) = df_i(\rho) \;\forall i \in I$.
\end{lem}

\begin{proof}
The differential of the group multiplication $f: G \times G \rightarrow G$ at
$(A_1,A_2) \in G \times G$ considered as a map $df(A_1,A_2): \mathfrak{g} \times
 \mathfrak{g}
\rightarrow \mathfrak{g}$ via right translation is given by
$$
df(A_1,A_2)(b_1,b_2)=b_1+Ad(A_1)b_2 \,.
$$
On the other hand the group multiplication law in $G \ltimes_{Ad} \mathfrak{g}$
reads
$$
(A_1,b_1)(A_2,b_2) = (A_1A_2, b_1 + Ad(A_1)b_2)\,,
$$
which proves the claim in this case. The general case is easily reduced to this
one. \hfill\qed
\end{proof}

\noindent In particular if the $f_i$ cut out $R(\Gamma,G)$ transversely,
then $T_{\rho}R(\Gamma,G)$ is identified with $Z^1(\Gamma;\mathfrak{g}_{Ad
  \circ \rho}$). In general, if $\rho\in R(\Gamma, G)$ is a smooth point, there
is a linear injection
\begin{equation}
\label{eqn:Weil}
\begin{array}{rcl}
T_{\rho}R(\Gamma,G) &\longrightarrow &Z^1(\Gamma;\mathfrak{g}_{Ad \circ \rho})\\

\dot{\rho}=\left.\textstyle\frac{d}{dt}\right\vert_{t=0}\rho_t &\longmapsto&
\dot{\rho}\rho^{-1}=\{\gamma
\mapsto \left.\textstyle\frac{d}{dt}\right\vert_{t=0} \rho_t(\gamma)\rho(\gamma)
^{-1}\}
\end{array}
\end{equation}
which identifies $T_{\rho}R(\Gamma,G)$ with a subspace of
$Z^1(\Gamma;\mathfrak{g}_{Ad \circ \rho})$.
This injection is known as \emph{Weil's construction}, and it extends naturally
to the non-smooth case by using the Zariski tangent space.
A cocycle $z$ is a coboundary, i.e.~$z \in B^1(\Gamma;\mathfrak{g}_{Ad
  \circ \rho})$, if and only if $z$ is tangent to the orbit of the action of
$G$ on $R(\Gamma,G)$ via conjugation.

\begin{lem}
\label{rep_var_smooth} $\rho_0 \in R(M,SU(2))$ is a smooth point
with tangent space $Z^1(\pi_1M;\R^3_{\rho_0})$. Similarly, $\rho_0
\in R(M,SL_2(\C))$ is a smooth point with tangent space
$Z^1(\pi_1M;\mathfrak{sl}_2(\C)_{Ad \circ \rho_0})$.
\end{lem}
\begin{proof}
There is an infinite sequence of obstructions for a cocycle $z \in
Z^1(\pi_1M;\R^3_{\rho_0})$ to be integrated into an actual path of
representations tangent to $z$, cf.~\cite{Gol1}. These obstructions
live in $H^2(\pi_1M; \R^3_{\rho_0})$. Using the injection of
Corollary~\ref{cor:injectiontoboundary} and the fact that the
obstructions are natural, we conclude that they vanish, because
they vanish in $H^2(\pi_1\partial\overline{M};\R^3_{\rho_0})$.
Then using Artin's theorem, cf.~\cite{Art}, they are in fact
integrable. Thus every element in $Z^1(\pi_1M;\R^3_{\rho_0})$ is
integrable, and therefore $\rho_0$ is a smooth point of
$R(M,SU(2))$ with tangent space $Z^1(\pi_1M;\R^3_{\rho_0})$.

For $SL_2(\mathbf C)$ the same argument applies since $\mathfrak{sl}_2(\C)_{Ad
  \circ \rho_0}$ is the complexification of the module $\mathfrak{su}(2)_{Ad
  \circ \rho_0} \cong \R^3_{\rho_0}$.
\hfill\qed
\end{proof}

\begin{prop}
$\chi_0 \in X(M,SU(2))$ is a smooth point of local real dimension
$q$ and tangent space isomorphic to $H^1(\pi_1M;\R^3_{\rho_0})$.
Similarly, $\chi_0 \in X(M,SL_2(\C))$ is a smooth point of local
complex dimension $q$ and tangent space isomorphic to
$H^1(\pi_1M;\mathfrak{sl}_2(\C)_{Ad \circ \rho_0})$.
\end{prop}

\begin{proof}
Since $\rho_0$ is irreducible, $Z(\rho_0(\pi_1M))=\{\pm 1\}$ in each
case. Since $SU(2)$ is compact (resp.~the action of $SL_2(\C)$ is proper on
the irreducible part of $R(M,SL_2(\C))$, cf.~\cite [Lemma 6.24]{Wei1}) and the
tangent space to the orbit is given by $Z^1(\pi_1M;\R^3_{\rho_0})$  (resp.~by
$Z^1(\pi_1M;\mathfrak{sl}_2(\C)_{Ad \circ \rho_0})$), the result follows with
Lemma \ref{rep_var_smooth}.\hfill\qed
\end{proof}

\noindent The fiber of the map
$$
ROT: R(M, \widetilde{\Isom \R^3}) \rightarrow R(M,SU(2))
$$
at $\rho \in R(M,SU(2))
$ is the the space of cocycles $Z^1(\pi_1M;\R^3_\rho)$, which is
given by a system of linear equations whose coefficients depend continuously on
$\rho$. Therefore its dimension is an upper semi-continuous function of $\rho$
in general.

If $\rho_0 \in R(M,SU(2))$ is a smooth point with $T_{\rho_0}
R(M,SU(2)) = Z^1(\pi_1M;\R^3_{\rho_0})$ (which is the case if the
$f_i$ cut out $R(M,SU(2))$ transversely or for $\rho_0 = ROT \circ hol$
according to Lemma \ref{rep_var_smooth}), then $T_\rho R(M,SU(2))$
injects into $Z^1(\pi_1M;\R^3_\rho)$ for $\rho$ in a neighbourhood of $\rho_0$
via Weil's construction (\ref{eqn:Weil}):
\begin{align*}
T_\rho R(M,SU(2)) &\longrightarrow Z^1(\pi_1M;\R^3_\rho)\\
\dot{\rho} &\longmapsto \dot{\rho}\rho^{-1}
\end{align*}
and hence $ROT: R(M, \widetilde{\Isom \R^3}) \rightarrow R(M,SU(2))$
is locally the projection of a vector bundle. More precisely we have:

\begin{lem}
$ROT: R(M,\widetilde{\Isom \R^3}) \rightarrow R(M,SU(2))$ is isomorphic to the tangent bundle $TR(M,SU(2))$ near $\rho_0$.
\end{lem}

\begin{proof}
The map
\begin{align*}
TR(M,SU(2)) &\longrightarrow R(M,\widetilde{\Isom \R^3})\\
(\rho,\dot{\rho}) &\longmapsto (\rho, \dot{\rho}\rho^{-1})
\end{align*}
is a vector bundle isomorphism near $\rho_0$.\hfill\qed
\end{proof}

\begin{cor}
$hol \in R(M,\widetilde{\Isom \R^3})$ is a smooth point.
\end{cor}

\noindent Let $\chi \in X(M,\widetilde{\Isom \R^3})$ denote
the character of $hol \in R(M,\widetilde{\Isom \R^3})$.

\begin{prop}\label{vectorbundle}
$\chi \in X(M,\widetilde{\Isom \R^3})$ is a smooth point of local real dimension
 $2q$.
The induced map
$$
ROT: X(M,\widetilde{\Isom \R^3}) \rightarrow X(M,SU(2))
$$
is locally the projection of a vector bundle isomorphic to the tangent bundle $T
X(M,SU(2))$ near $\rho_0$.
\end{prop}

\begin{proof}
Let us digress into a more general situation first:

Let $G$ be a compact Lie group, $M$ a smooth manifold (not
necessarily compact) and $G \times M \rightarrow M$ a smooth free
action. The associated infinitesimal action is the Lie algebra
homomorphism $\mathfrak{g}\rightarrow \Gamma(M,TM)$ defined by $ a
\mapsto \{ p \mapsto \left.\textstyle\frac{d}{dt}\right\vert_{t=0}
\exp(ta)p \}$. Since the original action was free, this
homomorphism is injective and we will identify $\mathfrak{g}$ with
its image in $\Gamma(M,TM)$.

Let $G\ltimes_{Ad} \mathfrak{g}$ be the semidirect product of $G$
and $\mathfrak g$ defined via $Ad: G \rightarrow
Aut(\mathfrak{g})$.  We extend the natural action of $G$ on $TM$
via the differential to an action of the group $G\ltimes_{Ad}
\mathfrak{g}$ in the following way:
\begin{align*}
(G \ltimes_{Ad} \mathfrak{g}) \times TM &\longrightarrow TM\\
((g,a),v) &\longmapsto dg(v) + a(g(\pi(v)))
\end{align*}
The relation $dg(a(p)) = (Ad(g)a)(gp)$ ensures that this defines a group
action. This action is free since the original action was free, it is clearly proper.

\begin{lem}\label{extendedaction}
$T(M/G)=TM/(G \ltimes_{Ad} \mathfrak{g})$.
\end{lem}

\begin{proof}
The fibers of the map $d\pi: TM \rightarrow T(M/G)$ are precisely
the orbits of the action of the group $G \ltimes_{Ad} \mathfrak{g}$ on $TM$. \hfill\qed
\end{proof}

\noindent We return to the situation of Proposition~\ref{vectorbundle}:

Clearly $SU(2) \ltimes_{Ad} \mathfrak{su}(2)\cong\widetilde{\Isom \R^3}$ via
$\mathfrak{su}(2)_{Ad} \cong \R^3$ as $SU(2)$-modules. For a representation $(\rho,z) \in
R(M,\widetilde{\Isom \R^3})$ and $(A,b) \in \widetilde{\Isom \R^3}$ we
have
$$
(A,b) (\rho,z) (A,b)^{-1} = (A\rho A^{-1}, Az + b - A\rho A^{-1}b) \,.
$$
Therefore the map
$
TR(M,SU(2)) \rightarrow R(M,\widetilde{\Isom \R^3})
$
is $\widetilde{\Isom \R^3}$-equivariant and Lemma \ref{extendedaction}
yields the proposition. \hfill\qed
\end{proof}

\noindent From the long exact cohomology sequence and Proposition \ref{vectorbundle} we obtain:

\begin{cor}\label{shortexact}
There is a short exact sequence
$$
0 \rightarrow H^1(M;TM) \rightarrow H^1(M;\E) \rightarrow H^1(M;TM)
\rightarrow 0\,.
$$
\end{cor}

\begin{rem}
Modulo the choice of a splitting of the exact sequence in
Corollary \ref{shortexact}, an infinitesimal deformation of the
holonomy of the Euclidean structure on $M$ is therefore determined
by the infinitesimal deformation of its rotational part and of its
translational part, and both can be independently prescribed.
\end{rem}

\noindent  Let $m_1,\ldots,m_q\in \pi_1M$ a system of meridians
for $\pi_1M$ (i.e.~one for each component of $Sing(C)$). We define
the angle function $\mu_j\!:\! U\subseteq X(M,SU(2))\to \mathbf
R$ in a neighborhood $U$ of $\chi_0$, so that
$\mu_j(\chi_0)=\alpha_j$. It is related to the trace by the
equality:
\[
\operatorname{trace}(\rho(m_j))=\pm
2\cos\frac{\mu_j(\chi_{\rho})}{2}\,.
\]
In particular $\mu_j$ is analytic. We extend $\mu_j$ to a neighborhood in
$X(M;SL_2(\mathbf  C))$ as a complex analytic function:
\[
\mu_j\!:\! V\subset X(M,SL_2(\mathbf  C))\to \mathbf  C\,,
\]
so that the complex length of $\rho(m_j)$ is $i\,\mu_j $ (i.e.~a
translation of length $Re(i\,\mu_j)$ plus a rotation of angle
$Im(i\,\mu_j)$).\\

\noindent The differentials $d\mu_j$ live in the cotangent space to
 the varieties of characters, thus
 \[
d\mu_j\in H^1(M;TM)^* \cong H_1(M;TM)\,.
 \]

\begin{prop}
\label{prop:kerdmu}
 $\ker \langle d\mu_1, \ldots , d\mu_q \rangle=\langle v
\rangle$, where $v$ is the class of the holonomy cocycle.
\end{prop}

\begin{proof}
To check that $d\mu_j(v)=0$, we observe that
$$H^1(\pi_1M;\R^3_{ROT \circ hol})\to
H^1(\langle m_j\rangle; \R^3_{ROT \circ hol})$$
maps $v$ to zero, because $hol(m_j)$
is a rotation with a fixed axis. In particular $d\mu_j$ evaluated
at $v$ is zero. To see the other inclusion, we use that
$$
\ker \langle d\mu_1, \ldots , d\mu_q \rangle = H^1(M;TM) \cap \left\langle [\omega^1_{len}], \ldots,
[\omega^q_{len}] \right\rangle \,,
$$
where we view $H^1(M;TM)$ as a subspace of $H^1(\partial\overline{M};TM)$ via
Lemma \ref{lem:inj}. Since the forms $\omega^j_{len}$ are $L^2$ near
the singular locus according to Lemma \ref{L2}, we conclude that $\ker\langle
d\mu_1,\ldots,d\mu_q\rangle\subseteq H^1_{L^2}(M;TM)$, which we view as a
subspace of $H^1(M;TM)$ via Corollary \ref{L2:inj}. Since $H^1_{L^2}(M;TM)$ is
spanned by $[id]$ the result follows. \hfill\qed
\end{proof}

\begin{cor}
\label{cor:dimdmu}
 \(\dim \langle d\mu_1,\ldots,d\mu_q\rangle=q-1. \)
\end{cor}

\section{Regeneration of structures}\label{section:regenerationstructures}
Along this section $\mathbf X^3$ denotes either $\mathbf S^3$ or $\mathbf H^3$.
The stabilizer of a point $p \in \mathbf X^3$ is isomorphic to
$SU(2)\cong Spin(
3)$.
Thus we view $\rho_0$ as a
representation of $\pi_1M$ in this stabilizer.

We consider a path of characters
\[
\begin{array}{rcl}
    [0,\varepsilon)&\to & X(M,\widetilde{\Isom\mathbf X^3})
    \\ t&\mapsto& \chi_t
\end{array}
\]
with $\chi_0=\chi_{\rho_0}$. Assume that this path is differentiable to the
right at $0$. The derivative $\frac{\partial\chi_t}{\partial t}(0)$ is
an element of $H^1(\pi_1M;\mathfrak g_{Ad \circ hol})$, where $\mathfrak g$
 is the Lie
algebra of $\widetilde{\Isom\mathbf X^3}$ (using Weil's
construction).

The Lie algebra $\mathfrak g$
decomposes into rotational and translational part with respect to
$p$:
\[
0\to \mathfrak t \to \mathfrak g \overset{trans}\longrightarrow T_p\mathbf
X^3\to 0
\]
where $\mathfrak t \cong \mathfrak{su}(2)$ is the Lie algebra of the stabilizer
of $p \in \mathbf X^3$. Thus
\[
trans \left(\frac{\partial \chi_t}{\partial t} (0)\right)\in
H^1(\pi_1M;\R^3_{ROT \circ hol})\,.
\]

\begin{thm}
\label{thm:regeneration}
 Let $\{\chi_t\}_{t\in [0,\varepsilon)}$ be a path
in $X(M,\widetilde{\Isom\mathbf X^3})$ with
$\chi_0=\chi_{\rho_0}$ and differentiable to the right at $0$.
If $trans(\frac{\partial\chi_t}{\partial t}(0))=
v$, then $\chi_t$ is the holonomy character of a $\mathbf X^3$-structure
on $M$, for $t\in (0,\delta)$ and some $0<\delta<\varepsilon$.

If in addition $\chi_t(m_j)$ is a rotation for each meridian
$m_j$, then the structure on $M$ completes to a cone manifold
structure.
\end{thm}

\begin{proof}
Consider $\mathbf X_t$, the space of constant sectional curvature
$\epsilon\, t^2$, where $\epsilon=\pm1$ is the curvature of $\mathbf X^3$,
for $t\in [0,\varepsilon)$.
Equivalently, $\mathbf X_t$ is $\mathbf X^3$ with the metric tensor scaled by
$t^{-2}$.
We fix a base point $p$ in $\mathbf X_t$
independently of $t$, so that the pointed Euclidean space $(\mathbf
X_0,p)$ is the limit of $(\mathbf X_t,p)$ when $t\to 0$, and consider the
union
$$
\overline {\mathbf X}=\operatornamewithlimits{\bigcup}\limits_{t\in
[0,\varepsilon]}\mathbf X_t.
$$
 We equip $\overline{\mathbf X}$ with a
manifold structure, using the local charts given by the exponential maps
at $p$ and the parameter $t\in [0,\varepsilon)$, after identifying
isometrically $T_p \mathbf X_t\cong\mathbf R^3$ for every $t\in
[0,\varepsilon)$. Notice that the choice of $p$ is relevant for
the topology  of $\overline{\mathbf X}$.

%For $t>0$ we have a diffeomorphism $h_t\! :\mathbf X_t\to
%\mathbf X_1=\mathbf X^3$ which is a dilation of factor $t$ centered at
%$p$:  $h_t(\exp_p(v))= \exp_p(t\, v)$ $\forall v\in T_p \mathbf
%X_t\cong \mathbf R^3\cong T_p\mathbf X_1$.

Choose a smooth path of representations $\rho_t$ with character $\chi_t$
so that $\rho_0=ROT \circ hol$ and $trans(\frac{\partial\rho_t}{\partial
t}(0))= trans(hol)$.

\begin{lem} For every $\gamma\in\pi_1 M$, the action of
$\rho_t(\gamma)$ on $\mathbf X_t$ for $t>0$ extends
continuously to the action of  $hol(\gamma)$ on $\mathbf X_0$ for the
$\mathcal C^1$-topology.
\end{lem}

Notice that $\rho_t$ acts on $\mathbf X_t$ isometrically, since
rescaling the metric does not change the isometry group.
The previous lemma provides an action on $\overline{\mathbf X}$,
which is isometric on each $\mathbf X_t$.

\begin{proof} We start describing the local coordinates. Let
%\begin{displaymath}
$
\exp_p^{(t)}\!:T_p\mathbf X_t\to\mathbf X_t
$
%\end{displaymath}
denote the Riemannian exponential.
Using the isometric identification
 $\mathbf R^3\cong T_p\mathbf X_t$, we have
\begin{displaymath}
   \exp_p^{(t)}(v)=\exp_p(t\, v)\qquad\forall  v\in \mathbf R^3,
\end{displaymath}
where $\exp_p=\exp_p^{(1)}:T_p\mathbf X^3\to \mathbf X^3$.
Thus, the (inverse of) the local charts is given as follows. Given an
open set $V\subset  \mathbf R^3$, $V\times
[0,\varepsilon)$ parametrizes a subset of $\overline{\mathbf X}$
via the map:
$$
\begin{array}{ccl}
    V\times \{0\} & \to & \mathbf X_0 \\
    (v,0) & \mapsto & v
\end{array}
\qquad
\begin{array}{ccl}
    V\times (0,\varepsilon) & \to &\overline{\mathbf X}\\
    (v,t) & \mapsto & \exp^{(t)}_p(v)=\exp_p(t\, v)\in\mathbf X_t
\end{array}.
$$
%%
%% , so that  $(v,0)\in V\times \{0\}$ is mapped to
%%$v\in V\subset \mathbf X_0$ and $(v,t)\in V\times (0,\varepsilon)$
%%to $\exp_p(t\, v)\in \mathbf X_t $.
To prove the lemma it suffices
to show that
$$
\lim_{t\to 0^+}
\frac1t\exp_p^{-1}(\rho_t(\gamma)\exp_p(t\,v))=hol(\gamma)(v)
$$ uniformly for $v\in  \mathbf R^3$ in a compact set for
the $\mathcal C^1$-topology.

The Lie algebra of $\operatorname{Isom}\mathbf
X_1$ decomposes as a   sum $\mathfrak g=\mathfrak t +\mathfrak
p$, where $\mathfrak t$ is the subalgebra of infinitesimal rotations
around $p$ and $\mathfrak p$ the subspace of infinitesimal translations
with respect to $p$. There is an isometric identification
$\mathfrak p\cong T_p\mathbf X_1\cong \mathbf R^3$, so that for any $v\in\mathbf
 R^3$:
$$
\exp_p(v)=\exp(v)(p),
$$
here we view $v\in T_p\mathbf X_1\cong  \mathbf R^3$ when
we write $\exp_p(v)$ and $v\in\mathfrak p$ for $\exp(v)(p)$.

%% we view $v\in\mathbf R^3$ as an element of both $ T_p\mathbf
%%X_1$  and $   \mathfrak p$.

According to global Cartan's decomposition, we write
$\rho_t(\gamma)=\exp (b_t)\, a_t$, where $b_t\in \mathfrak p$
and $a_t$ belongs to the stabilizer of $p$ in $\mathbf X_1$.
Notice that $\frac{b_t}t\to trans(hol(\gamma)) $ as $t\to 0$ by hypothesis.

Using this notation:
\begin{multline*}
\rho_t(\gamma)(\exp_p(t\,v))=(\exp(b_t)\, a_t\, \exp(t\,  v)) (p)=
(\exp(b_t)\, a_t\, \exp(t\,  v)\, a_t^{-1}) (p)=\\
\exp (b_t)\exp(t\, a_t(v))(p)
=\exp (b_t + t\,  a_t (v) + t^2\, C)(p)
\end{multline*}
where $C=C(b_t/t,a_t,v, t)$ is an  analytic function. Here we use that $b_t=O(t)
$
together with the Campbell-Hausdorff formula.
Using again this formula and the fact that $a_t(v), b_t\in\mathfrak p$, we get:
\begin{displaymath}
\rho_t(\gamma)(\exp_p(t\,v))=\exp_p (b_t + t\, a_t (v) + t^2\,  C'),
\end{displaymath}
where $C'=C'(b_t/t,a_t,v, t)$ is also  analytic.
Hence
$$
\frac1t\exp_p^{-1}(\rho_t(\gamma)\exp_p(t\, v))=a_t(v)+\frac{b_t}t+ t \, C'(b_t/
t,a_t,v, t) .
$$
When $t\to 0$ this converges to $hol(\gamma)(v)$ uniformly for $v$
in a compact subset of $\mathbf R^3$, because $ \frac{b_t}t\to trans(hol(\gamma)
)$ and
$a_0=\rho_0(\gamma)$ is a lift of $ROT(hol(\gamma))$, by
hypothesis.

To prove convergence in the $\mathcal C^1$-topology, we write $v+\varepsilon \,
w$
with $v,w\in\mathbf R^3$, and $\varepsilon>0$ small. The previous calculation
yields easily:
\begin{displaymath}
        \frac1t\exp_p^{-1}(\rho_t(\gamma)\exp_p(t\,(v+\varepsilon \, w)))=
        \\ a_t(v+\varepsilon \, w)+\frac{b_t}t+
        t\, C'(b_t/t,a_t,v+\varepsilon \, w, t)\,.
\end{displaymath}
We can compute the derivative with respect to $\varepsilon$
when $\varepsilon=0$ and we get $\mathcal C^1$-convergence.
\hfill\qed
\end{proof}

\noindent We modify slightly Goldman's construction to deform the structure.
Let $\overline M$ be a compact core of $M=C^{smooth}$. Consider
$E=\widetilde {\overline M}\times_{\pi_1M} \overline{\mathbf X}$,
which is a bundle over $\overline M$ with fiber $\overline{\mathbf
X}$. It is in fact a union of bundles $E_t=\widetilde {\overline
M}\times_{\pi_1M}  {\mathbf X}_t$ with fiber ${\mathbf X}_t$. Each
$E_t$ has a natural flat connection that varies continuously with
$t$. The developing map  of the Euclidean structure induces a
section $s:M\to E$ with values in $\mathbf X_0$ transverse to the flat
connection. Since the image of $s$ is compact, we compose it with
the flow $\Phi_t$ of the vector field tangent to the direction of
$t$ (this is defined globally when $\epsilon=-1$ but not when
$\epsilon =1$, because $\mathbf X_0$ is not homeomorphic to
$\mathbf X_t$). For small values of $t>0$, $\Phi_t\circ
s\!:\overline M\to E_t $ is a section, still transverse to the
flat connection by $\mathcal C^1$-continuity, hence inducing a developing map.

This provides a structure on $\overline M$, that can be completed
by controlling its behavior on $\partial\overline M$, using the
hypothesis about the meridians. \hfill\qed
\end{proof}

\noindent The volume of these cone manifolds is increasing with $t$, because the
 starting
Euclidean structure is viewed as totally degenerate with volume zero. In fact we
 get more precisely:

\begin{prop}[Schl\"afli's formula]
    \label{prop:Schlafli}
        Let $C_{t}$ be the family of cone manifolds constructed in
        Theorem~\ref{thm:regeneration} of constant curvature $K\in\mathbf  R$.
        Assume moreover that the path of characters is analytic.
        Then
        \[
        {K}\,d \operatorname{vol}{C_t}=\frac12\sum_e
        \operatorname{lenght} (e)d\alpha_e
        \]
where the sum runs over the edges and circles $e$ of $Sing(C)$.
\end{prop}

\noindent\emph{Sketch of proof.} The construction of developing
maps in the proof of Theorem~\ref{thm:regeneration} can be made so
that we have a set of points $Z=\{z_1\ldots,z_k\}$ such that the
balls $B_{r_i}(z_i)$ cover $C$, where $r_i$ is much smaller that
the injectivity radius, and $\mathcal D_t(\tilde z_i)$ varies
analytically (we assume that the path of characters is analytic).
Let $C_t$ the cone structure with holonomy $\chi_t$ and  define
\[
P_i(t)=\{x\in C_t\mid d(x,z_i)\leq (x,z_j)\ \forall j=1,\ldots k\}
\]
By construction, this $P_i$ is an analytic family of polyhedra in
the space $\mathbf X^3$. Analyticity implies that the topological
type of $P_i$ changes in a discrete subset of times $t$. Thus we
may apply Schl\"afli's formula to them. Adding all the terms, we
get the formula of the proposition, because the contribution of
nonsingular edges is trivial. See \cite[Prop.~4.2]{Porti} for
further details. \hfill\qed

\begin{rem}
\label{rem:trans} In the hyperbolic case, $\widetilde{\Isom\mathbf
H^3}\cong SL_2(\mathbf  C)$ and $\mathfrak{sl}_2(\mathbf
C)=\mathfrak{su}(2)\otimes_{\mathbf  R}\mathbf  C$. Therefore in
this case $trans:H^1(\pi_1 M; \mathfrak g_{Ad \circ hol}) \to
H^1(\pi_1M; \R^3_{ROT \circ hol})$ can be viewed as the imaginary
part.

In the spherical case, $\widetilde{\Isom\mathbf S^3}\cong
SU(2)\times SU(2)$ and the Lie algebra of the stabilizer of a
point $p \in \mathbf S^3$ is conjugate to the diagonal subalgebra
of $\mathfrak{su}(2)\times \mathfrak{su}(2)$. Therefore in this
case $trans:H^1(\pi_1 M;\mathfrak g_{Ad \circ hol}) \to
H^1(\pi_1M; \R^3_{ROT \circ hol})$ can be viewed as the difference
between factors.
\end{rem}

\begin{rem}
\label{rem:rotation}
For every $\rho$ with character $\chi_{\rho}\in V\subset X(M,SL_2(\mathbf
C))$, \[ \rho(m_i)\textrm{ is a rotation if and only if }
\mu_i(\chi_\rho)\in\mathbf  R. \]

For every $(\rho_1,\rho_2)$ with character
$(\chi_{\rho_1},\chi_{\rho_2}) \in  U\times U \subset
X(M,SU(2))\times X(M,SU(2))$, \[ (\rho_1(m_i),\rho_2(m_i))\textrm{
is a rotation if and only if }
\mu_i(\chi_{\rho_1})=\mu_i(\chi_{\rho_2}).
\]
\end{rem}

\section{Constructing paths of hyperbolic and spherical structures}\label{sec:pathsofstructures}
Up to changing the indices, by Corollary~\ref{cor:dimdmu} we may
assume that $d\mu_2,\ldots,d\mu_q$ are linearly independent. Thus the
set
\[
\mathcal C=\{\chi\in U\subset X(M,SL_2(\mathbf  C))\mid
\mu_i(\chi)=\mu_i(\chi_0), \ i\geq 2\}
\]
is a smooth complex curve in a neighborhood $U$ of $\chi_0$.

Notice that the class of the holonomy cocycle $v$ is the tangent
vector to $\mathcal C$ at $\chi_0$, by
Proposition~\ref{prop:kerdmu}. Thus, using Remark~\ref{rem:trans}
we have:

\begin{rem} To every path $\gamma:[0,\varepsilon)\to \mathcal C$
with $\gamma(0)=\chi_0$ and $\operatorname{Im}(\gamma'(0))\neq 0$,
Theorem~\ref{thm:regeneration} applies.
\end{rem}

\begin{lem}
\label{lem:deg2}
The restricted map
$\mu_1\vert_{\mathcal C}\! :\!\mathcal C\to\mathbf C$ is a branched covering of
degree two.
\end{lem}

\begin{proof}
By Proposition~\ref{prop:kerdmu}, $d\mu_1\vert_{\mathcal
C}(\chi_0)=0$, thus $\mu_1\vert_{\mathcal C}$ is either constant or a
branched covering of degree $d\geq 2$. Seeking a contradiction,
assume that $\mu_1\vert_{\mathcal C}$ is constant, then we have a
path of characters in the curve $\mathcal C $ to which we apply
Theorem~\ref{thm:regeneration}. Since each $\mu_i$ stays constant on
this curve, this path corresponds to hyperbolic cone manifolds with
constant cone angle. Thus by Schl\"afli's formula
(Proposition~\ref{prop:Schlafli}) they have constant volume,
contradicting the fact that the volume increases from zero. Thus
$\mu_1\vert_{\mathcal C}$ is a branched covering of degree $d\geq 2 $.

To prove that the degree $d$ is precisely $ 2$,  we assume that $d>2 $ and
seek again  a contradiction. We look at the inverse image or real points
$(\mu_1\vert_{\mathcal C})^{-1}(\mathbf  R)$, because for characters here,
the image of $\mu_1$ is a rotation. This inverse image is a graph with $2
d$ branches starting at $\chi_0$, and the angle between the branches is
$\pi/d$. Thus two of the branches are real and the remaining $2d-2$ have
nontrivial imaginary part. Hence Theorem~\ref{thm:regeneration} applies to
those $2d-2$ branches. In addition, if $d>2$, then there are branches for
which $\mu_1$ is strictly larger than $\mu_1(\chi_0)$ and branches for
which $\mu_1$ is strictly less than $\mu_1(\chi_0)$. Since $\mu_1$ is the
cone angle, we have constructed regenerating families of hyperbolic cone
manifolds with both increasing and decreasing  cone angles, contradicting
Schl\"afli's formula. \hfill\qed
\end{proof}

\begin{cor}
\label{cor:hyperbolicregeneration}
There is a family of hyperbolic cone structures obtained by
decreasing $\alpha_1$.
\end{cor}

\begin{proof}
Since $\mu_1\vert_{\mathcal C}$ is a branched covering of degree 2,
$(\mu_1\vert_{\mathcal C})^{-1}(\mathbf  R)$ has 4 branches starting at
$\chi_0$. Two of them are real, and to the other two one applies
Theorem~\ref{thm:regeneration}. The cone angle $\alpha_1$ of those
branches must decrease by Schl\"afli's formula. These two branches correspond to
complex conjugate representations, i.e.~with opposite orientations,
because changing the sign of $v$ corresponds to changing the orientation. \hfill
\qed
\end{proof}

\begin{prop}
There is a family of spherical cone structures obtained by
increasing $\alpha_1$.
\end{prop}

\begin{proof}
To construct the spherical structure we consider
\[
\mathcal D=\mathcal C\cap X(M,SU(2)).
\]
This is a real analytic curve in a neighborhood $V$ of $\chi_0$.
It is the set of real points of $\mathcal C$.

By Lemma~\ref{lem:deg2}, the map $\mu_1\vert_{\mathcal D}\! :\!\mathcal
D\to\mathbf  R$ is locally equivalent to the map $x\mapsto x^2 $ in a
neighborhood of $0\in\mathbf  R$. Thus the inverse fiber of the map
$\mu_1\vert_{\mathcal D}$ consists of pairs of points, except for
$\mu_1(\chi_0)$, which consists of a single point. Let $\mathcal
D^{\pm}\subset \mathcal D$ be two subintervals of $\mathcal D$ such that
$\mathcal D^+\cup \mathcal D^-=\mathcal D$ and $\mathcal D^+\cap \mathcal
D^-=\{\chi_0\}$. We consider the set
\[
S=\{(\chi_+,\chi_-)\in \mathcal D^+\times\mathcal D^- \mid
\mu_1(\chi_+) =\mu_1(\chi_-)\}
\]
%%
%%
%% Thus $S$ has
%%precisely four branches starting at $\chi_0$, two of them with
%% $\chi_+=\chi_-$ and two other with $\chi_+\neq\chi_-$.
%%Theorem~\ref{thm:regeneration} applies to the latter two branches,
%%because the difference of tangent vectors gives the translational
%%part, that is necessarily a non-zero multiple of $v$ (see
%%Remark~\ref{rem:trans}).
%%
%%
%%Now consider the set
%%\[
%%S=\{(\chi_+,\chi_-)\in \mathcal D\times\mathcal D  \mid
%%\mu_1(\chi_+) =\mu_1(\chi_-) \}
%%\]
because $\mu_1^+(\chi_+)=\mu_1^-(\chi_-)$ is the condition that
guarantees that $(\chi_+,\chi_-)$ is the character of a
representation $\rho$ in $Spin(4)$ so that $\rho(\mu_1)$ is a
rotation, by Remark~\ref{rem:rotation}.
%%
%%
%%
%%Then we consider the map
%%\[
%%(\mu_1^+,\mu_1^-)\! :\!\mathcal D\times\mathcal D\to\mathbf
%%R\times \mathbf  R
%%\]
%%where $\mu_1^{\pm}$ denotes $\mu_1\vert_{\mathcal D}$ on each
%%component of the product $\mathcal D\times\mathcal D$. Now
%%consider the set
%%\[
%%S=\{(\chi_+,\chi_-)\in \mathcal D\times\mathcal D  \mid
%%\mu_1^+(\chi_+)=\mu_1^-(\chi_-) \}
%%\]
%%because $\mu_1^+(\chi_+)=\mu_1^-(\chi_-)$ is the condition that
%%guaranties that $(\chi_+,\chi_-)$ is the character of a
%%representation $\rho$ in $Spin(4)$ so that $\rho(\mu_1)$ is a
%%rotation, by Remark~\ref{rem:rotation}.
%%
%%By Lemma~\ref{lem:deg2}, the map $\mu_1\vert_{\mathcal D}\!
%%:\!\mathcal D\to\mathbf  R$ is locally equivalent to the map
%%$x\mapsto x^2 $ in a neighborhood of $0\in\mathbf  R$. Thus $S$ has
%%precisely four branches starting at $\chi_0$, two of them with
%% $\chi_+=\chi_-$ and two other with $\chi_+\neq\chi_-$.
Theorem~\ref{thm:regeneration} applies to the latter two branches,
because the difference of tangent vectors gives the translational
part, that is necessarily a non-zero multiple of $v$ (see
Remark~\ref{rem:trans}). \hfill\qed
\end{proof}

\noindent If we replace $\mathcal D^+\times \mathcal D^-$ by $\mathcal
D^-\times \mathcal D^+$ (i.e.~we change the order of the factors)
then we get the same structures with different orientation.

\begin{rem}
\label{rem:furtherregenerations}
 We have constructed regenerations
by deforming the cone angle $\alpha_1$ because we assumed that the
forms
\[
d\mu_2,\ldots,d\mu_q
\]
are linearly independent. Below in
Corollary~\ref{cor:linearrelation} we prove that the linear
relation satisfied by these forms is:
\[
\sum_{i=1}^q l_i \, d\mu_i=0\qquad\text{ for some } l_i>0.
\]
Thus we can regenerate by deforming any of the cone angles.
\end{rem}

\section{Fold locus and Euclidean structures}
\label{sec:foldlocus}
Now we analyze the behavior of the real analytic map
$$
\mu=(\mu_1,\ldots,\mu_q) :
X(M,SU(2))\textsc{}\to\mathbf  R^q
$$
in a neighborhood of $\chi_0$.

\begin{lem}
\label{lem:equivalencex2}
In a neighborhood of $\chi_0$, $\mu$ is equivalent to
the map
\[
(x_1,x_2,\ldots,x_q)\mapsto (x_1^2,x_2,\ldots,x_q)
\]
in a neighborhood of the origin via a real analytic change of coordinates at both the source and at the target.
\end{lem}

\begin{proof}
By Corollary~\ref{cor:dimdmu} and up to permuting coordinates, we
may assume that the forms $d\mu_2,\ldots,d\mu_q\in T^*_{\chi_0}
X(M,SU(2))$ are linearly independent. Thus the
set
\[
\mathcal {D}=\{\chi\in U\subset X(M,SU(2))\mid
\mu_i(\chi)=\mu_i(\chi_0), \ i\geq 2\}
\]
is a smooth (real analytic) curve in a neighborhood $U$ of $\chi_0$.
By Lemma~\ref{lem:deg2} we can choose real analytic
coordinates $(x_1, \ldots, x_q)$ in a neighborhood of
$\chi_0$ with $x_i(\chi_0)=0$ such that
$
\mu_1\vert_{\mathcal{D}} = \mu_1(\chi_0) + x_1^2
$
and
$
\mu_i=\mu_i(\chi_0) + x_i
$
for $i \geq 2$. Note that the curve $\mathcal{D}$ corresponds to the $x_1$-axis
in these coordinates.

Since $\partial \mu_1/\partial x_1(0)=0$ and $\partial^2
  \mu_1/\partial x_1^2(0) = 1$, the implicit function theorem yields that the
  set $\{ \partial \mu_1/\partial x_1 = 0 \}$ is locally around $0$
  the graph of a real analytic function $f=f(x_2, \ldots, x_q)$. Let $g=g(x_2,
  \ldots, x_q)$ be the real analytic function defined by
$$
g(x_2, \ldots, x_q) = \mu_1(f(x_2,
  \ldots,x_q), x_2, \ldots, x_q)\,.
$$
Let further $h$ be the unique real analytic
  function satisfying
$$h^2(x_1, \ldots,
  x_q)=\mu_1(x_1, \ldots, x_q) - g(x_2,\ldots, x_q)
$$
and
$
\partial h/\partial x_1(0)=1
$.
Then the map $\phi$ defined by
$$
\phi(x_1, \ldots, x_q) = (h(x_1, \ldots , x_q), x_2, \ldots, x_q)
$$
is a local diffeomorphism.
We obtain $\mu_1 \circ \phi^{-1} = x_1^2 + g (x_2,\ldots,x_q)$
and $\mu_i \circ \phi^{-1} = x_i+\mu_i(\chi_0)$ for $i \geq 2$. With the local
  diffeomorphism $\psi$ defined by
$$
\psi(x_1, \ldots, x_q) = (x_1-g(x_2, \ldots, x_q), x_2-\mu_2(\chi_0), \ldots, x_q-\mu_q(\chi_0))
$$
we obtain $\psi \circ \mu_1 \circ \phi^{-1}(x_1, \ldots, x_q) =
(x_1^2,x_2,\ldots,x_q)$.
%x_1^2$ and $\psi \circ \mu_i \circ \phi^{-1}(x_1, \ldots, x_q) =
%x_i$ for $i \geq 2$.
\hfill\qed
\end{proof}

\begin{rem}
The same result holds true for the complex analytic extension
$$
\mu=(\mu_1,\ldots,\mu_q) :
X(M,SL_2(\C))\textsc{}\to\mathbf  C^q
$$
just by composing with the complex analytic extensions of the coordinate changes.
\end{rem}

\begin{dfn}
The \emph{fold locus} $\Fold\subset X(M,SU(2))$ is  the set of
points where $\mu$ is not a local diffeomorphism.
\end{dfn}

\noindent In the coordinate system of Lemma~\ref{lem:equivalencex2} the
fold locus is the set
 $\{x_1=0\}$, hence:

\begin{cor}
\label{cor:Jmanifold} There exists  a neighborhood $U$ of $\chi_0$
such that both $\Fold\cap U$ and $\mu(\Fold\cap U)$ are  codimension one
submanifolds. In addition, for each $\chi\in \Fold$, $\dim \ker\langle
d\mu_1,\ldots,d\mu_k\rangle=1$.
\end{cor}

\begin{prop} Each $\chi\in \Fold\cap U$ is the rotational part of the
holonomy of a Euclidean cone manifold structure. The translational part lies
in $\ker\langle d\mu_1,\ldots,d\mu_k\rangle$.
\end{prop}

\begin{proof} When we deform a character $\chi$ in $\Fold$ we can also
deform continuously the vector $v$ in $\ker\langle
d\mu_1,\ldots,d\mu_k\rangle\subset T_{\chi}X(M,SU(2))$, which is a
one dimensional subspace that varies continuously, by
Lemma~\ref{lem:equivalencex2}.
 This corresponds to deforming continuously the holonomy $hol$ by a
family of representations of $\pi_1M$ in
$\widetilde{\Isom\R^3}$. The condition $v\in\ker
\langle d\mu_1,\ldots,d\mu_k\rangle$ ensures that the image of meridians
are rotations (cf. \cite[Prop.~9.6]{Porti}).
 Since those representations
map the meridians to rotations, it follows that they correspond to
holonomy representations of Euclidean cone manifolds. \hfill\qed
\end{proof}

\begin{cor}
\label{cor:Eeucl}
Points of $E=\mu(\Fold)$ are angles of Euclidean cone manifold structures on $C$.
\end{cor}

\noindent Let $\bar l=(l_1,\dots,l_q)$ denote the lengths of the singular circles
and components of $C$. Those are unique up to homothety.

\begin{prop}
\label{prop:lnormal}
The vector $\bar l$ is normal to $E=\mu(\Fold)$ at $\mu (\chi_0)$.
\end{prop}

\begin{proof}
This is a consequence of Schl\"afli's formula
(Prop.~\ref{prop:Schlafli}). \hfill\qed
\end{proof}
\begin{cor}
\label{cor:linearrelation} On $T_{\chi_0}(X(M,SU(2)))$,
$\sum_{i}l_i\,d\mu_i=0$.
\end{cor}

\begin{proof}
It holds on $T_{\chi_0}\Fold$ by Proposition~\ref{prop:lnormal}. In
addition, it also holds on $v=[trans\circ hol]$. By
Lemma~\ref{lem:equivalencex2} and Proposition~\ref{prop:kerdmu},
those spaces generate the whole $T_{\chi_0}(X(M,SU(2))$. \hfill\qed
\end{proof}

\noindent It follows from this corollary that we can regenerate into
hyperbolic or spherical cone manifold structures by decreasing or
increasing any of the
cone angles, the geometry  depends on the sign of $\sum
l_i\alpha'_i$. (See Remark~\ref{rem:furtherregenerations}.)

\begin{cor}[Local Rigidity for Euclidean cone manifolds]
\label{cor:LRECM}
 Let $C$ be a closed Euclidean cone manifold with cone angles $\leq \pi$ which is
not almost product.
 Then
deformations up to dilations of $C$ into Euclidean cone manifold
structures
 are parametrized by the $q$ cone
angles in a manifold $E$ of dimension $q-1$ and transverse to the
vector of singular lengths $\bar{l}=(l_1,\ldots,l_q)$.
\end{cor}

\begin{figure}[ht]
\psfrag{F}{$\mathcal{F}$}
\psfrag{pi}{$\pi$}
\psfrag{0}{$0$}
\psfrag{E}{$E$}
\psfrag{mu}{$\mu$}
\psfrag{v}{$v$}
\psfrag{X(M,SU(2))}{$X(M,SU(2))$}
{\centering\
        \epsfig{file =  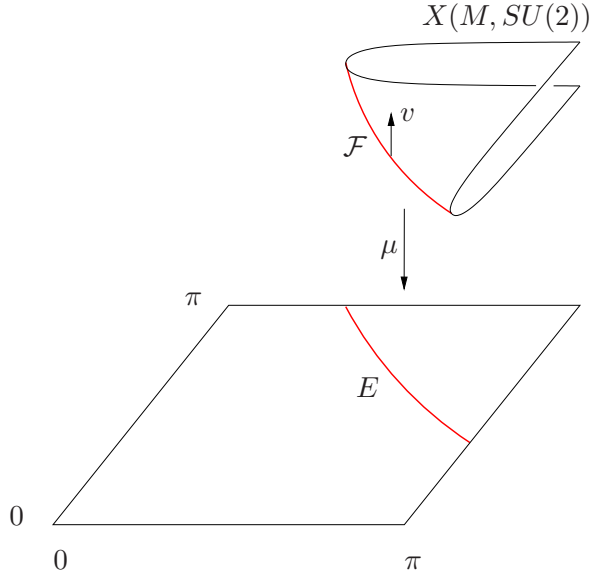, width = .5\textwidth}\caption{
 Geometry of the map $\mu$ and the fold locus.}\label{fig:foldlocus}
\par}
\end{figure}

\begin{proof}
 The previous analysis determines all representations of
$\pi_1M$ into $\widetilde{\Isom \R^3}$ up to conjugation that map
meridians to rotations, i.e.~pairs $(v,\chi_{\rho})$ such that
$\chi_\rho\in X(M,SU(2))$ and further $v\in \ker \langle
d\mu_1,\ldots,d\mu_k\rangle\subset T_{\chi} X(M,SU(2))$, see Figure \ref{fig:foldlocus}. \hfill\qed
\end{proof}

\begin{cor}[Local Rigidity for regenerations]
\label{cor:LRR} Let $C$ be as above. Then deformations of $C$ into
constant curvature structures up to dilations are  parametrized by
the $q$ cone angles, and the type of geometry depends on the side
of $E$.
\end{cor}

\begin{proof}
The analysis of this section and previous ones determine the
structure of the spaces of characters of $\pi_1M$ in
$\widetilde{\Isom \mathbf X^3}$. Namely, spherical regenerations
correspond to pairs of $SU(2)$-characters $(\chi_+,\chi_-)$ such
that $\mu(\chi_+) = \mu(\chi_- )$. According to Lemma
\ref{lem:equivalencex2}, these pairs are determined by the value
of $\mu$, which is the multiangle of the structure, cf.~Figure
\ref{fig:foldlocus}. Hyperbolic regenerations correspond to
$SL_2(\C)$-characters $\chi$ with $\mu(\chi) \in \R^q$, which are
not $SU(2)$-characters, cf.~Figure \ref{fig:hypreg}. Again such a
character is determined up to complex conjugation by the value of
$\mu$.

\begin{figure}[ht]
\psfrag{F}{$\mathcal{F}$}
\psfrag{v}{$v$}
\psfrag{iv}{$iv$}
\psfrag{X(M,SU(2))}{$X(M,SU(2))$}
\psfrag{hypreg}{hyperbolic regenerations (pos.~or.)}
{\centering\
         \epsfig{file =  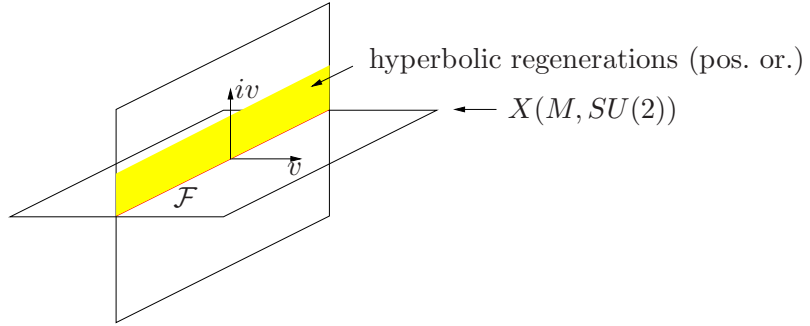, width = .5\textwidth}
         \caption{The locus of $SL_2(\C)$-characters having real traces.}\label{fig:hypreg}
\par}
\end{figure}

\noindent Next, we shall show that the deformations of
the structures are uniquely determined by the deformations of the characters,
adapting an argument from
\cite{Gol2}.

Let $\overline M$ be a compact core of $M=C^{smooth}$, so that we
will consider deformations of structures on $\overline M$, knowing
that the structure on the surface $\partial\overline M$ is
determined by the holonomy, and the structure on $M-\overline
M$ is given by taking the cone.

Using the notation of
Section~\ref{section:regenerationstructures}, $E=\cup_{t\in [0,t)}
E_t$ is the union of flat bundles $E_t$ with fiber $\mathbf X_t$
the space of constant curvature $\pm t^2$, and holonomy the
deformation of the Euclidean holonomy representation.  Since $E_t$
is flat, there is a horizontal foliation, which is also used to
define locally a projection onto the fibers $\mathbf X_t$.

The Euclidean structure is induced by a section  $s_0\! :
\overline M\to E_0$ transverse to the horizontal foliation. In
Section~\ref{section:regenerationstructures}, $s_0$ is deformed to
$s_t\! :\overline M\to E_t$, following the flow in the time
direction on $E$. Let $\sigma_t\! :\overline M\to E_t$ be another
deformation, i.e.~a section of class $\mathcal C^1$ such that
$\sigma_0=s_0$ and $\sigma_t\vert_{\partial\overline M}=s_t
\vert_{\partial\overline M}$. We claim that, for $t\geq 0$
sufficiently small, there exist a diffeomorphism $h_t$ of $M$
isotopic to the identity so that the developing maps corresponding
to $\sigma_t$ and $s_t\circ h_t$ are equivalent (namely, the
sections $\sigma_t$ and  $s_t\circ h_t$ composed with the local
projections to $\mathbf X_t$ differ by isometries of $\mathbf
X_t$). Following  \cite{Gol2}, we take a foliated neighborhood $N$
of $s_0(\overline M)$ in $E_0$, so that the intersection of leaves
with $N$ are balls, and consider $N\times [0,t)$ in $E$, via the
flow as in Section~\ref{section:regenerationstructures}. For $t$
sufficiently small, closeness and transversality imply that the
image of both $s_t$ and $\sigma_t$ meet each horizontal leave in
$N \times \{t\}$ precisely once. The intersection with the
horizontal leaves can be used to define the map $h_t$ as in
\cite{Gol2}. \hfill\qed
\end{proof}

\section{Global results}
\label{sec:global}

In this section we finish the proof of Theorem~\ref{thm:main}. We first
establish some properties of the set of Euclidean angles and recall
some results about the set of hyperbolic and spherical ones.

Let $E$ denote the set of multiangles with  a Euclidean cone manifold structure.

\begin{prop}
\label{prop:completesubmfld}
 $E$ is a properly embedded hypersurface of
$(0,\pi)^q$.
%%(Here completeness means that  $E$ can be enlarged so
%%that $(\overline E-E)\cap (0,\pi)^q=\emptyset$.)
\end{prop}

\begin{proof}
First we show that every Euclidean cone manifold structure on $C$
with cone angles in $(0,\pi)$ is not almost product. This will
imply that $E$ is  a $q-1$ submanifold, by
Corollary~\ref{cor:Eeucl}.
By contradiction, we assume that there exists such a Euclidean cone structure
on $C$ which is almost product.
Since all cone angles are in $(0,\pi)$,
this implies that $C$ is Seifert
fibered and $Sing(C)$ is a union of fibers. By
Corollary~\ref{cor:almostseifert}, in this case $C$ is almost
product for any cone angle, contradicting the hypothesis of Theorem~\ref{thm:main}.

To prove properness, we show that if $x_n$ is a sequence in $E$
converging to $x_{\infty}\in (0,\pi)^{q}$, then $x_{\infty}\in E$.
Let $C_n$ be the corresponding sequence of cone structures on $C$
with multiangles $x_n$. We can rescale them so that
$\operatorname{diam} C_n=1$. Since the cone angles are uniformly
$<\pi$, by Theorem 5.3 and Corollary 5.4 of \cite{BLP2} the cone manifolds $C_n$
 are
uniformly thick: they have a base point $x_n\in C_n$ with
$inj(x_n)\geq c>0$. Thus they have a convergent subsequence to a
Euclidean cone manifold of diameter 1, hence a cone structure on $C$
with cone angles $x_{\infty}$. \hfill\qed
\end{proof}

\noindent In \cite[Ch.~3]{Wei2} the following proposition is proved, which
applies to $C$ under the assumptions of the paper.

\begin{prop}
 \label{prop:largerspherical}
 If $\bar \alpha\in (0,\pi]^q$ is a
multiangle of a spherical cone structure on $C$ which is not
Seifert fibered, then every $\bar
\beta=(\beta_1,\dots,\beta_q)\in (0,\pi]^q$ with
$\alpha_i\leq\beta_i$ is the multiangle of a spherical cone
structure with cone angles $\bar\beta$.

If $\bar \alpha\in (0,\pi]^q$ is a multiangle of a hyperbolic cone
structure on $C$, then every $\bar
\beta=(\beta_1,\dots,\beta_q)\in (0,\pi]^q$ with
$0<\beta_i\leq\alpha_i$ is the multiangle of a hyperbolic cone
structure with cone angles $\bar\beta$.
\end{prop}

\begin{rem}
\label{rem:smallness}
 When all cone angles are $\pi$, the orbifold is spherical and hence small
(without essential 2-suborbifolds). In particular all turnovers in $C$ are compressible
or boundary parallel.
\end{rem}

\noindent
Before proving Theorem~\ref{thm:main}, we need to establish incompatibility
between structures of different sign. We normalize the constant curvature to be
$-1$, $0$ or $1$.

\begin{lem}
\label{lem:incompatibility}
Given $\bar\alpha=(\alpha_1,\ldots,\alpha_q)\in(0,\pi]^q$, $\bar\alpha$
cannot be the multiangle of two cone structures with different
normalized constant curvatures  ($-1$, $0$ or $1$).
\end{lem}

\begin{proof}
Seeking a contradiction, assume for instance that $\bar\alpha\in(0,\pi]^q$ is the multiangle of a hyperbolic
and Euclidean cone structure simultaneously.
Using Corollary~\ref{cor:hyperbolicregeneration}, the Euclidean cone structure can be regenerated
to hyperbolic ones, which, by construction,
have arbitrarily small diameter, and also cone multiangle arbitrarily
close to $\bar\alpha$. On the other hand, when we perturb
the hyperbolic cone structure with cone multiangle $\bar\alpha$, we find hyperbolic
cone manifolds with the same multiangles but diameter bounded below away from zero. This contradicts
the global rigidity of
hyperbolic cone manifolds proved in \cite[Thm. 1.4]{Wei2}.
The same argument works for spherical
and Euclidean (global rigidity in the spherical case is the content of \cite[Thm.\ 1.7]{Wei2}).

Finally we deal with the case that $\bar\alpha$ is simultaneously
the multiangle of a hyperbolic and spherical structure on $C$.
Consider $n\in\mathbf N$ sufficiently large so that $2\pi/n<\alpha_i$ for $i=1,\ldots,q$. By
Proposition~\ref{prop:largerspherical},
$(\frac{2\pi}n,\ldots,\frac{2\pi}n)$ is the multiangle of hyperbolic cone manifold, hence a hyperbolic orbifold.
Notice that this orbifold cannot be spherical, because its fundamental group is
infinite. Consider now a path of cone multiangles
between $\bar\alpha$ and $(\frac{2\pi}n,\ldots,\frac{2\pi}n)$  which is decreasing along each component.
 By Proposition~\ref{prop:largerspherical} every multiangle in this path corresponds to a hyperbolic cone structure.
 However, sphericity
must fail at some
multiangle of this path, because
$(\frac{2\pi}n,\ldots,\frac{2\pi}n)$  is not the multiangle of a spherical structure.
By the arguments in the proof of
Proposition~\ref{prop:completesubmfld}, when sphericity fails we find precisely
a Euclidean structure, and we have reduced to the first case of incompatibility
between hyperbolic and Euclidean
structures. \hfill\qed
\end{proof}

\noindent
The following proposition is also used in the proof of Theorem~\ref{thm:main}.

\begin{prop}
\label{prop:angle0}
Every multiangle in $[0,\pi)^q$ with at least one angle $0$ is the
multiangle of a hyperbolic cone structure.
\end{prop}

\begin{proof}
By \cite[Cor.~1.5]{Wei2}, the smooth part
of $C$ is hyperbolic.
Applying Thurston's Dehn
filling theorem, we have the proposition for a neighborhood of the
origin. To cover the rest of the multiangles, we just have to enlarge
some of the angles, keeping the other ones equal to zero (i.e.~complete cusps).
Thus we have a lower bound on the cone angles,
and in addition the diameter is infinity. Applying the results of
\cite{BLP2}, we can  enlarge each one of the cone angles up to $<\pi$.
Notice that some convergence results of \cite{BLP2} apply to cone manifolds 
without essential turnovers.
In our case there are no turnovers by Remark~\ref{rem:smallness} except if we allow
cusped turnovers (i.e.\
turnovers with some cone angle 0). However cusped turnovers are not a problem in
 those arguments, because
all cone angles are $<\pi$, and therefore they cannot converge to a Euclidean cone 2-manifold.
\hfill\qed
\end{proof}

\begin{rem}
\label{rem:opennessat0}
The hyperbolic structures of the previous proposition can be deformed in a neighborhood in
$[0,\pi)^q$.
\end{rem}

\noindent
This remark can be easily proved adapting the arguments of the proof of hyperbolic Dehn
filling for orbifolds in \cite[App. B]{BP} and using the infinitesimal rigidity
results established in
\cite{Wei1}.

\medskip

\noindent\emph{Proof of Theorem~\ref{thm:main}:}
Assume first that $C$ is a Euclidean cone manifold as in the statement and all
cone angles are $\pi$. Then using the regeneration results of Section~\ref{sec:pathsofstructures}
and Proposition~\ref{prop:largerspherical}, every point  in $(0,\pi)^q$ is the multiangle
of a hyperbolic cone structure on $C$.

Assuming that  at least one of the angles of $C$ is $<\pi$, then there exists a
spherical cone structure
on $C$ with all cone angles $\pi$, again by Section~\ref{sec:pathsofstructures}
and Proposition~\ref{prop:largerspherical}.
We consider all segments  in $[0,\pi)^q$ starting at $(\pi,\ldots,\pi)$ and ending at some point with
at least one coordinate zero.
The first point of the segment is the multiangle of a spherical cone structure and the last one
is the multiangle of a hyperbolic one, by Proposition~\ref{prop:angle0}.
By  Remark~\ref{rem:opennessat0}, we can assume that the multiangle
of the hyperbolic structure
lies in the open cube $(0,\pi)^q$.
Starting at $(\pi,\cdots,\pi)$, we move along the segment by decreasing the cone
 angles and obtaining a family
of spherical cone manifolds. This family cannot be spherical all the time because the endpoint  of the
segment corresponds to a hyperbolic structure and they are incompatible by Lemma
~\ref{lem:incompatibility}.
 Since all cone angles are $<\pi$, then the end of the spherical
subsegment is the multiangle of a Euclidean structure, by
the same argument as in Proposition~\ref{prop:completesubmfld}.
By the regeneration results of Section~\ref{sec:pathsofstructures}
and Proposition~\ref{prop:largerspherical},
we connect the multiangle of this Euclidean structure to the endpoint
of the segment by multiangles of hyperbolic cone structures.

The previous argument shows that every point in  $(0,\pi)^q$ is the multiangle of
a constant curvature cone structure on $C$, and that $E$ is a hypersurface that
divides
$(0,\pi)^q$ in two components $H$ and $S$ corresponding respectively to hyperbolic and
 spherical multiangles.

For the uniqueness, notice first that in Lemma~\ref{lem:incompatibility} we establish incompatibility between hyperbolic,
spherical and Euclidean for a given multiangle.
Global rigidity for hyperbolic and spherical structures is proved in \cite{Wei2}. Using local rigidity of the regenerations
(Corollary~\ref{cor:LRR}) and global rigidity for the hyperbolic
structures, we get global rigidity for the Euclidean ones.
\hfill \qed

\section{The Whitehead link}
\label{sec:example}

\begin{figure}[h]
{\centering\
        \epsfig{file =  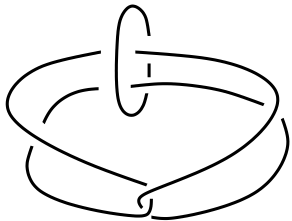, width = .25\textwidth}\caption{The Whitehead link.}\label{fig:whitehead}
\par}
\end{figure}

In this section we illustrate the main theorem for cone manifold structures on the 3-sphere with singular locus given by the Whitehead link, which is the 2-component link depicted in Figure \ref{fig:whitehead}. Furthermore we discuss geometric structures corresponding to multiangles contained in the boundary of $[0,\pi
]^2$ to some extent.

\noindent R.N.~Shmatkov has computed the curve of multiangles of
Euclidean structures by constructing the corresponding cone
manifolds explicitly, cf.~\cite{Shm}. These structures are not
almost product.  The hyperbolic region had earlier been computed
by A.~Mednykh in a similar way, cf.~\cite{Med} and references
therein.

Alternatively, one can proceed as follows: Let $M$ denote the
complement of the Whitehead link in $S^3$.  If $a, b \in \pi_1M$
are meridians around the two components of $\Sigma$, the
fundamental group has the following presentation, cf.~\cite{HLM}:
$$
\pi_1M = \langle a,b \,\vert\, awa^{-1}w^{-1} \rangle
$$
with $w=bab^{-1}a^{-1}b^{-1}ab$.
The $SL_2(\mathbf C)$-character variety of $M$ has been computed in \cite{HLM}.
Namely, after identifying $X(M,SL_2(\mathbf C))$ with the image of the map
\begin{align*}
(t_a,t_b,t_{ab}): R(M,SL_2(\mathbf C)) &\rightarrow \C^3\\
\rho & \mapsto (\operatorname{tr} \rho(a), \operatorname{tr} \rho(b),\operatorname{tr} \rho(ab))
\end{align*}
in $\mathbf C^3$, it is given
by
$$
X(M,SL_2(\mathbf C))=\{ (x,y,z) \in \mathbf C^3 : p(x,y,z) \cdot q(x,y,z)=0\}
$$
with
$$
p(x,y,z) = xy - (x^2+y^2-2)z+xyz^2 - z^3
$$
and
$$
q(x,y,z)=x^2+y^2+z^2-xyz-4\,.
$$

Again by \cite{HLM}, the irreducible part of the character
variety, i.e.~those  characters which correspond to irreducible
representations, is given by
$$
X^{irr}(M,SL_2(\mathbf C)) = \{ p(x,y,z)=0\} \setminus \{q(x,y,z)=0\}\,;
$$
moreover
$$
\{ p(x,y,z) = 0\} \cap \{ q(x,y,z) = 0)\} = \{x=\pm 2, z=\pm y\} \cup \{y=\pm 2,
 z=\pm x\}\,.
$$

Holonomy representations of hyperbolic cone manifold structures lift to irreducible
 representations $\rho$ with $t_a(\rho) \in (-2,2)$ and $t_b(\rho) \in (-2,2)
$. We may write
$$
x = \pm 2 \cos (\alpha/2)
$$
and
$$
y = \pm 2 \cos(\beta/2)
$$
with $\alpha$ and $\beta$ the cone angles around the two components of $\Sigma$.

\begin{figure}[ht]
\psfrag{(0,p)}{$(0,\pi)$}
\psfrag{(p2,p)}{$(\pi/2,\pi)$}
\psfrag{(p,p2)}{$(\pi,\pi/2)$}
\psfrag{(p,p)}{\textcolor{blue}{$(\pi,\pi)$}}
\psfrag{(p,0)}{$(\pi,0)$}
\psfrag{(0,0)}{$(0,0)$}
\psfrag{hyperbolic}{ hyperbolic}
\psfrag{euc}[][][1][-45]{\textcolor{red}{Euclidean}}
\psfrag{sph.}{\textcolor{blue}{spherical}}
{\centering\
        \epsfig{file =  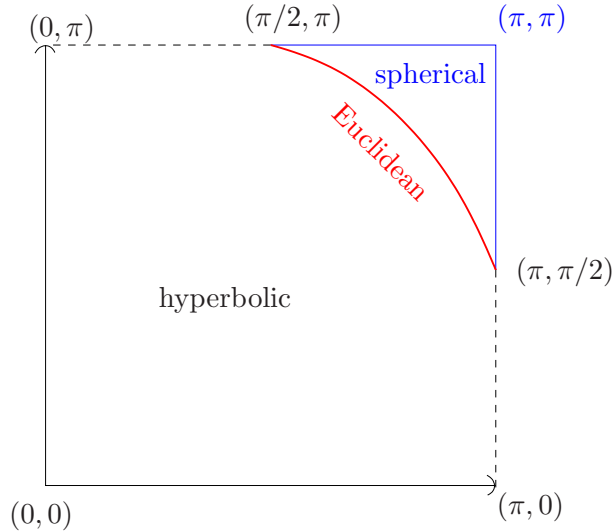, width = .5\textwidth}\caption{
 The curve of Euclidean multiangles and the hyperbolic and spherical regions.}
 \label{fig:euclideanregion}
\par}
\end{figure}

Rotational parts of Euclidean holonomies correspond to
representations as above where in addition $t_a$ and $t_b$ fail to
be local coordinates. This is precisely the locus where the
discriminant of $p$ computed with respect to the variable $z$
vanishes. The discriminant is given by the following polynomial:
\begin{align*}
f(x,y)=x^6y^2-2x^4y^4+2x^4y^2+x^2y^6+2x^2y^4-11x^2y^2+32\\
-48x^2-48y^2+24y^4+24x^4-4x^6-4y^6 \,.
\end{align*}
For the resulting curve of multiangles of Euclidean cone manifold structures see
 Figure \ref{fig:euclideanregion}.

For multiangles contained in some parts of the boundary of
$[0,\pi]^2$ we can construct cone manifold structures as well
using our main theorem: For multiangles of the form $(0,\beta)$,
$0 \leq \beta < \pi$, resp.~$(\alpha,0)$, $0\leq \alpha < \pi$, we
obtain hyperbolic cone manifold structures, whereas for
multiangles
 of the form  $(\pi,\beta)$, $\pi/2<\beta\leq \pi$, resp.~$(\alpha,\pi)$, $\pi/2
 < \alpha \leq \pi$, we obtain spherical ones.

A $Nil$-orbifold structure with branching indices $(4,2)$,
resp.~$(2,4)$, i.e.~corresponding to multiangles $(\pi/2, \pi)$,
resp.~$(\pi,\pi/2)$, has been constructed by E.~Su\'arez,
cf.~\cite{Sua}.

For the remaining part of the boundary we do not give a complete description, we
 rather prove the following statement:

\begin{lem}
There is no hyperbolic cone manifold structure corresponding to
multiangles contained in $[0,\pi/2) \times \{\pi\} \cup \{\pi\} \times [0,\pi/2)$.
\end{lem}

\begin{proof}
Let us suppose there existed such a structure. Then we consider
the double branched cover of $S^3$ branched along the component of
the Whitehead link with cone angle $\pi$. This branched cover is
again $S^3$ since the components of the Whitehead link are
unknotted. The other component of the Whitehead link lifts to the
 torus link $T(4,2)$, whose complement is known to be non-hyperbolic. On the
 other hand the lift of the intital hyperbolic cone metric will be a hyperbolic cone
 metric on $S^3$ with singular locus $T(4,2)$, which is a contradiction in view
of the results in \cite{Wei2}. \hfill\qed
\end{proof}

\begin{rem}
Multiangles of the form $(2\pi/n,\pi)$, resp.~$(\pi, 2\pi/n)$,
with $n\geq5$ correspond to orbifold structures modelled on
$\widetilde{PSL_2(\mathbf R)}$-geometry.
\end{rem}

%%\section{Further comments}
%%
%%
%%
%%
%%When some of the cone angles are $\pi$, there are other
%%possibilities. The Spherical do not degenerate, but the hyperbolic
%%and Euclidean ones can.
%%
%%
%%There are other examples e.g. Borromean rings.

\sc{
Departament de Matem\`atiques, Universitat Aut\`onoma de
Barcelona, \\ E-08193 Be\-lla\-ter\-ra, Spain, porti@mat.uab.es}

\sc{
Mathematisches Institut, Universit\"at M\"unchen, Theresienstra\ss
e 39, \\ D-80333 M\"un\-chen, Germany,
weiss@mathematik.uni-muenchen.de
}

\end{document}